\documentstyle[12pt,epsf]{article}

\begin{document}

\newcommand{\rhob}{\overline{\rho}}
\newcommand{\ud}{{\rm d}}
\newcommand{\ui}{{\rm i}}
\newcommand{\ue}{{\rm e}}

\title{\bf On the spacing distribution of the Riemann zeros: corrections to the
asymptotic result}

\author{\it E. Bogomolny\dag, O. Bohigas\dag, P. Leboeuf\dag \ and A. G.
Monastra\ddag \\ {\small\dag Laboratoire de Physique Th\'eorique et Mod\`eles
Statistiques \thanks{Unit\'e de recherche associ\'ee au CNRS.} , Universit\'e
de Paris XI,} \\ {\small B\^at. 100, 91405 Orsay Cedex, France}\\ {\small
\ddag TU Dresden Institut f\"ur Theoretische Physik, 01062 Dresden, Germany}}


\maketitle

\begin{abstract}
It has been conjectured that the statistical properties of zeros of the
Riemann zeta function near $z = 1/2 + \ui E$ tend, as $E \rightarrow \infty$,
to the distribution of eigenvalues of large random matrices from the Unitary
Ensemble. At finite $E$ numerical results show that the nearest-neighbour
spacing distribution presents deviations with respect to the conjectured
asymptotic form. We give here arguments indicating that to leading order these
deviations are the same as those of unitary random matrices of finite
dimension $N_{\rm eff}=\log(E/2\pi)/\sqrt{12 \Lambda}$, where $\Lambda=1.57314 \ldots$ is
a well defined constant.
\end{abstract}

\section{Introduction}

\begin{figure}
\begin{center}
\leavevmode
\epsfxsize=6.7in
\epsfysize=2.5in
\epsfbox{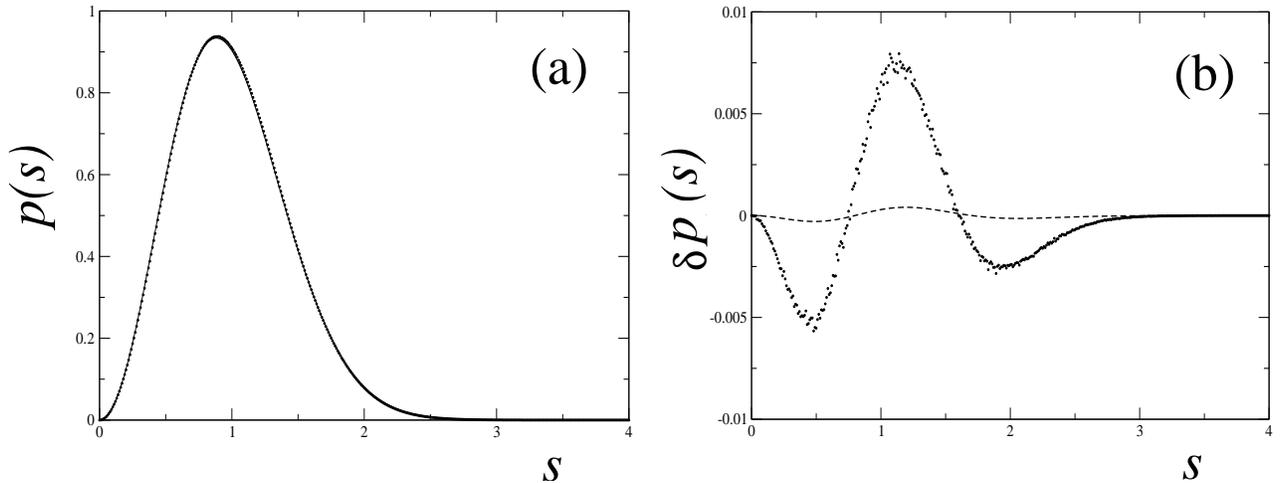}

\caption{{\small Nearest neighbour spacing distribution of the Riemann zeros
located in a window near $E=2.5041178 \times 10^{15}$. (a) Numerical results
(dots, from A. Odlyzko) compared to the asymptotic CUE curve (full line),
almost indistinguishable. (b) Difference between the numerical result and the
asymptotic CUE curve (dots) compared to the difference between the spacing
distribution of CUE matrices of size $N_{0}$ (Eq.~(\ref{N0})) and the
asymptotic curve (dashed line).}}

\end{center}
\label{Fig1}
\end{figure}

\noindent 
The study of connections between random matrix theory and properties of the
Riemann zeta function, $\zeta(z)$, has known recently significant developments
\cite{jpa}. A central point is Hugh Montgomery's (generalized) conjecture
\cite{m} that in the asymptotic limit (high on the critical line for $z=1/2 +
\ui E$) the fluctuation properties of non-trivial Riemann zeros are the same
as for the circular unitary ensemble (${\rm CUE_N}$) of $N \times N$ random
unitary matrices (with Haar measure) in the limit of large dimensionality,
$N\to \infty$. In particular, the normalized pair correlation function of the
Riemann zeros with $E\to \infty$ is conjectured to be
\begin{equation} \label{r2gue}
R_2 (s) = 1 -  \left( \frac{\sin  \pi s}{\pi s}  \right)^2\; ,
\end{equation}
where $s$ is the unfolded distance between zeros (i.e. the mean spacing is set
to one). Andrew Odlyzko, since the late 70's, started accurate and extensive
numerical computations of Riemann zeros in order to check this and other
conjectures \cite{od1}. His main result is that in the limit $E\to \infty$
correlation functions of Riemann zeros do agree with random matrix
predictions. For instance, in Fig.~\ref{Fig1}(a) the density distribution
$p(s)$ of spacings among consecutive zeros of $\zeta(z)$ (called the
nearest--neighbour spacing distribution in random matrix literature) is
plotted for a billion zeros around the $10^{16}$-th zero. The agreement with
the CUE asymptotic prediction $p_0 (s)$ (represented by a solid line in that
figure) is remarkable.

It is interesting to look at the difference $\delta p (s)=p(s)-p_0(s)$ between
the computed and the conjectured distributions (see Fig.~\ref{Fig1}(b)).
Though the difference is small (of order $10^{-2}$) it has a clear structure
with a nontrivial $s$ dependence.

One may wonder how this difference compares with the one obtained, within
random matrix theory, between the asymptotic and a finite $N$ calculation. To
proceed, one needs a criterion for the size of the matrix to compare with. The
simplest assumption is to choose $N=N_0$, where
\begin{equation}\label{N0}
N_{\scriptscriptstyle {0}} = \log \left( \frac{E}{2\pi} \right) \ .
\end{equation}
This matrix size is obtained by equating the local density of zeros at height
$E$ along the critical line to the density of eigenvalues of the $N\times N$
unitary matrix \cite{ks}.

This dimensional correspondence has been successful when comparing statistical
properties of $\zeta (1/2 + \ui E)$ at finite $E$ with those of characteristic
polynomials of ${\rm CUE_N}$ matrices of size $N=N_0$ \cite{ks}. The Riemann
zeros in Fig.~\ref{Fig1} are located in a window around $E=2.5041178 \times
10^{15}$, which gives $N_{0} = 33.6188$. The difference between the finite $N$
(with $N=N_0$) and the asymptotic nearest neighbour spacing distributions is
represented by a dashed line in Fig.~\ref{Fig1}(b). Though the functional form
of the correction is qualitatively correct, its amplitude is clearly too small
(by a factor of order 20). In his paper \cite{od1}, A. Odlyzko commented:
'Clearly there is structure in the difference graph, and the challenge is to
understand where it comes from'. The purpose of this note is to provide some
elements in this direction.

\section{Two--point correlation functions}

Let us first analyze the corrections to the asymptotic two--point function
(\ref{r2gue}) for $N\times N$ random matrices from the circular unitary
ensemble ${\rm (CUE_N)}$ and for zeros of the Riemann zeta function. The
comparison between them will provide an effective set of rules or
correspondences between both sequences.

We start first with random matrix theory. For $N$-dimensional unitary matrices
the $n$--point correlation functions are given by \cite{r2}
\begin{equation} \label{rn}
R_n^{\scriptscriptstyle {\rm (CUE_N)}} (x_1,\ldots ,x_n)=\det (K(x_i,x_j))|_{i,j=1\ldots ,n}\ ;
\end{equation}
$K(x_i,x_j)$ is the corresponding kernel,
\begin{equation} \label{kern}
K(x,y)=\frac{\sin(\pi(x-y))}{N\sin(\pi
(x-y)/N)}=K_0(x-y)+\frac{1}{N^2}K_1(x-y) +{\cal O} ( N^{-4} )
\end{equation}
where
\begin{equation}
K_0(s)= \frac{\sin(\pi s)}{\pi s}\ 
\end{equation}
and 
\begin{equation}
K_1(s)= \frac{\pi s}{6} \sin(\pi s)\ . 
\end{equation}  
In particular, the unfolded two--point correlation function takes the form
\begin{equation} \label{cuen}
R_2^{\scriptscriptstyle {\rm (CUE_N)}} (s)=
1 - \left ( \frac{\sin (\pi s) }{N\sin ( \pi s/N)}\right )^2 \ ;
\end{equation}
the expansion (\ref{kern}) leads to
\begin{equation} \label{r2cue}
R_2^{\scriptscriptstyle {\rm (CUE_N)}} (s) = 1 - \frac{\sin^2 (\pi s) }{ \pi^2
s^2 } - \frac{ 1 }{ 3 N^2 } \sin^2 (\pi s) - \frac{(\pi s)^2}{ N^4 } \sin^2
(\pi s) + {\cal O} ( N^{-6} )\ .
\end{equation}
This formula expresses the correlation function as the asymptotic result
Eq.~(\ref{r2gue}) plus corrections proportional to inverse even powers of the
matrix dimension. Notice that the corresponding expansion for random hermitian
matrices (Gaussian Unitary Ensemble) is different, it is non--analytic and
includes a term proportional to $N^{-1}$ (cf \cite{r2,tw}).

We now need an equivalent result for the Riemann zeros. An heuristic formula
for the two-point correlation function for these zeros was obtained by
Bogomolny and Keating in Ref.\cite{bk} using the Hardy-Littlewood conjecture
for the distribution of prime pairs (for more details see
\cite{Varenna,LesHouches}). It states that the two-point correlation function
of Riemann zeros, $r_2 (\epsilon)$, is the sum of three terms
\begin{equation}
r_2 (\epsilon) = \rhob^2 + r_2^{\rm (diag) } (\epsilon) + r_2^{\rm (off)
} (\epsilon) \ ,
\end{equation}
where the smooth density of zeros, $\rhob$, is asymptotically  
\begin{equation}
\rhob=\frac{1}{2\pi}\log \left (\frac{E}{2\pi}\right )\ ,
\end{equation}
and the diagonal, $r_2^{\rm (diag) } (\epsilon)$, and off-diagonal, 
$r_2^{\rm (off)} (\epsilon) $, parts are given by the following convergent
expressions 
\begin{equation}\label{diag}
r_2^{\rm (diag) } (\epsilon) = -\frac{1}{4 \pi^2} \frac{\partial^2}{\partial
\epsilon^2} \left[ \log | \zeta( 1 + \ui \epsilon ) |^2 + 2 \sum_p
\sum_{r=1}^{\infty} \frac{ 1- r }{ r^2 p^{r} } \cos(\epsilon r \log p) \right]
\end{equation}
and
\begin{equation}\label{off}
r_2^{\rm (off) } (\epsilon) = \frac{1}{4 \pi^2} | \zeta( 1 + \ui
\epsilon ) |^2 \ue^{\ui 2 \pi \rhob \epsilon } \prod_p \left[ 1
-\frac{ (1 - p^{\ui \epsilon} )^2 }{ (p - 1)^2 } \right] + {\rm c. c.} \ .
\end{equation}
Here the summation and the product are taken over all primes $p$. In
\cite{LesHouches} it was checked that these formulas agree very well with
Odlyzko's results for the two-point correlation function of Riemann zeros.

The unfolded two--point correlation function is obtained by measuring
distances between zeros in units of the local mean spacing,
\begin{equation} \label{r2u}
R_2 (s) = \frac{1}{\rhob^2 } r_2 \left (\frac{s}{\rhob}\right ) \ .
\end{equation}
We are interested in the corrections to the asymptotic behavior of $R_2 (s)$
in the limit when $E \rightarrow \infty$. In this limit $\rhob \rightarrow
\infty$ and the argument of $r_2$ in Eq.(\ref{r2u}) becomes small (keeping $s$
finite). Therefore, one can expand $r_2 (\epsilon)$ for $\epsilon \ll 1$. To
perform the expansion it is convenient to use the known series
$$
\zeta( 1 + x) = \frac{1}{x} +\sum_{n=0}^{\infty}\frac{(-1)^n}{n!}\gamma_n
x^n\ ,
$$
where $\gamma_i$ are the Stieljes constants, and the following auxiliary
expansions:
$$ 
\frac{\partial^2}{\partial \epsilon^2} \sum_p \sum_{r=1}^{\infty}
\frac{ 1- r }{ r^2 p^{r} } \cos(\epsilon r \log p) =
\sum_{n=0}^{\infty} c_n \epsilon^{2 n}
$$
where
$$ 
c_n = \frac{(-1)^n}{ (2 n)! } \sum_p (\log p)^{2(n+1)}
\sum_{r=1}^{\infty} \frac{(r-1) r^{2 n} }{p^r} \ , 
$$
and 
$$
\prod_p \left[ 1 -\frac{ (1 - p^{\ui \epsilon} )^2 }{ (p - 1)^2 }
\right] = 1 + c_0 \epsilon^2 + \ui Q \ \epsilon^3 + {\cal O}
(\epsilon^4) 
$$
with 
$$ Q = \sum_p \frac{ \log^3 p }{ (p-1)^2 } \ .$$

Collecting the different terms in the expansion, one gets
\begin{equation}
r_2^{\rm (diag) } (\epsilon) = - \frac{1}{2 \pi^2 \epsilon^2 } -
\frac{(\gamma_0^2 + 2 \gamma_1 + c_0) }{2 \pi^2 } + {\cal O}
(\epsilon^2)
\end{equation}
and
\begin{equation}
r_2^{\rm (off) } (\epsilon) = \frac{1}{4 \pi^2} \left[ \frac{1}{\epsilon^2 } +
(\gamma_0^2 + 2 \gamma_1 + c_0) + \ui Q \epsilon + {\cal O} (\epsilon^2)
\right] \ue^{\ui 2 \pi \rhob \epsilon } + {\rm c. c.}\ .
\end{equation}
The unfolded two--point correlation function takes therefore the form
\begin{equation} \label{r2r1}
R_2 (s) = 1 - \frac{\sin^2 (\pi s) }{ \pi^2 s^2 } - \frac{ (\gamma_0^2 + 2
\gamma_1 + c_0) }{\pi^2 \rhob^2} \sin^2 (\pi s) - \frac{ Q }{2 \pi^2
\rhob^3} s \sin (2 \pi s) + {\cal O} ( \rhob^{-4} )\ .
\end{equation}
Equation (\ref{r2r1}) expresses the two--point correlation function of Riemann
zeros as the asymptotic random matrix result given by Eq.~(\ref{r2gue}), plus
corrections that are proportional to inverse powers of the average density of
zeros.

The comparison of Eqs.(\ref{r2cue}) and (\ref{r2r1}) shows that the leading
terms coincide, as conjectured by Montgomery. To relate sub-leading terms we
proceed as follows. Up to ${\cal O} ( \rhob^{-4} )$ in Eq.~(\ref{r2r1}), the
term of ${\cal O} (\rhob^{-3})$ can be absorbed in the term of ${\cal O} (
\rhob^{-2} )$ by rescaling the $s$ variable. Equation (\ref{r2r1}) now takes
the form
\begin{equation} \label{r2r2}
R_2 (s) = 1 - \frac{\sin^2 (\pi s) }{ \pi^2 s^2 } - \frac{ (\gamma_0^2 + 2
\gamma_1 + c_0) }{\pi^2 \rhob^2} \sin^2 (\pi \alpha s) + {\cal O} (
\rhob^{-4}) \ ,
\end{equation}
where
\begin{equation} \label{alpha}
\alpha = 1 + \frac{Q}{2 \pi \rhob \ (\gamma_0^2 + 2 \gamma_1 + c_0)}
=1+\frac{C}{\log (E/2\pi)} \ ,
\end{equation}
$C = Q/(\gamma_0^2 + 2 \gamma_1 + c_0) = 1.4720...$. The comparison of
Eqs.~(\ref{r2cue}) and (\ref{r2r2}) leads to the following conclusions:
\begin{itemize}
\item[i)] To leading order in $1/\rhob$ the two--point correlation function of
the Riemann zeros coincides with that of eigenvalues of random ${\rm CUE_N}$
matrices of effective dimension $N=N_{\rm eff}$, where
\begin{equation}\label{neff}
N_{\rm eff} = \frac{\pi \rhob }{\sqrt{3 \Lambda} } = \frac{1}{\sqrt{12 \Lambda}} \log
\left( \frac{E}{2 \pi} \right) \ ,
\end{equation}
where $\Lambda \equiv \gamma_0^2 + 2 \gamma_1 + c_0 = 1.57314 \ldots$. 
\item[ii)] The next-to-leading order is obtained by rescaling the variable $s$
in the first correction term according to
\begin{equation}\label{sca}
s \rightarrow \alpha s \ ,
\end{equation}
where $\alpha$ is given by Eq.(\ref{alpha}).
\end{itemize}
The effective matrix size obtained from our analysis is, therefore, different
from $N_{0}$ (Eq.~(\ref{N0})), the relation between them being a
multiplicative factor $N_{\rm eff} = (12 \Lambda)^{-1/2} N_{0} = 0.230158
N_{0}$.

\section{The nearest--neighbour spacing distribution}

In general, to compute the nearest--neighbour spacing distribution it is
necessary to know correlation functions of arbitrary number of points. Though
for the Riemann zeros the Hardy-Littlewood conjecture is sufficient to obtain
all correlation functions \cite{bk2}, the computations are cumbersome and for
simplicity we shall use another method.

For CUE, the nearest--neighbour spacing distribution may be expressed as
\cite{r2}
\begin{equation} \label{psdef}
p^{\scriptscriptstyle {\rm (CUE_N)}} (s)=\frac{d^2 E(s)}{ds^2}
\end{equation}
where
\begin{equation} \label{enofs}
E (s) = \det \left[ \delta_{j k} - \frac{\sin \left( \frac{\pi
s}{N} (j - k) \right)}{\pi (j - k)} \right] \ , 1 \leq j, k \leq N \ . 
\end{equation}

Some physical arguments \cite{aa} indicate that to leading order, under quite
general conditions deviations from standard random matrix theory in systems
with no time reversal symmetry are reduced to a change of the kernel only, and
that the correction term is the same as for ${\rm CUE_N}$ matrices with an
effective matrix size.

\begin{figure}
\begin{center}
\leavevmode
\epsfxsize=6.6in
\epsfysize=2.4in
\epsfbox{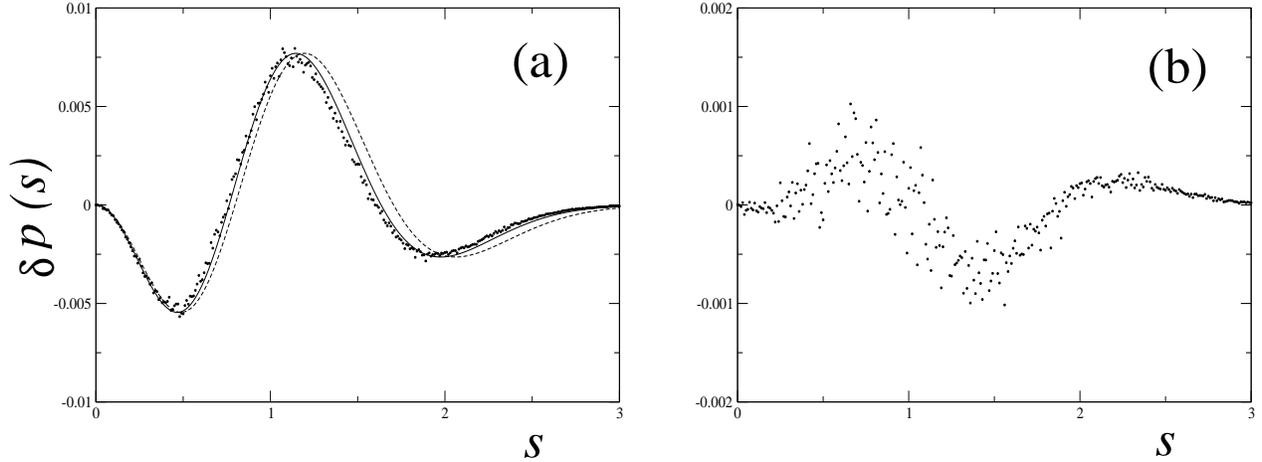}

\caption{{\small (a) Difference between the nearest neighbour spacing
distribution of the Riemann zeros and the asymptotic CUE distribution for a
billion zeros located in a window near $E=2.5041178 \times 10^{15}$ (dots),
compared to the theoretical prediction Eq.(\ref{psz2}) (full line). The dashed
line does not include the scaling of $s$. (b) Difference between the numerical
values for the Riemann zeros (dots) and the full curve (theory) of part (a).}}

\end{center}
\label{Figg2}
\end{figure}

If we accept this conjecture then from the results of the previous section it
follows that to leading order {\it all} correlation functions of the Riemann
zeros are the same as those of ${\rm CUE_N}$ matrices with effective dimension
given by Eq.~(\ref{neff}). In particular, it means that the nearest--neighbour
spacing distribution of the Riemann zeros can be calculated as follows. First,
find the expansion of the nearest--neighbour distribution for ${\rm CUE_N}$
random unitary matrices in inverse powers of $N$, namely
\begin{equation}\label{pscue}
p^{\scriptscriptstyle {\rm (CUE_N)}} (s) = p_0 (s) + \frac{1}{N^2} \ p_1 (s) +
{\cal O} ( N^{-4} ) \ .
\end{equation}
(The fact that the first correction is of order $N^{-2}$ follows from results
in \cite{tw}). The expansion (\ref{pscue}) is difficult to derive
analytically. We have therefore computed numerically the spacing distribution
$p^{\scriptscriptstyle {\rm (CUE_N)}} (s)$ from
Eqs.~(\ref{psdef}-\ref{enofs}). The correction term is then obtained by
computing $p_1 (s) = N^2 \ [p^{\scriptscriptstyle {\rm (CUE_N)}} (s) - p_0
(s)]$ for increasing values of $N$. An alternative method would be to use a
non-linear differential equation for $p^{\scriptscriptstyle {\rm (CUE_N)}}
(s)$ as derived in \cite{tw}.

Second, to leading order, at a given height $E$ above the real axis, replace
in Eq.~(\ref{pscue}) $N$ by $N_{\rm eff}$ given by Eq.~(\ref{neff}). Finally,
we shall approximate the next-to-leading correction of $p (s)$ according to
the rule (\ref{sca}). In such an approximation the nearest--neighbour spacing
distribution of the Riemann zeros equals the universal random matrix result,
$p_0(s)$, plus the correction $\delta p(s)$, where
\begin{equation} \label{psz2}
\delta p (s) = \frac{1}{N_{\rm eff}^2} \ p_1 (\alpha s) + {\cal O} ( N_{\rm eff}^{-4} )\ .
\end{equation}
$N_{\rm eff}$ and $\alpha$ are given by Eq.~(\ref{neff}) and (\ref{alpha}),
respectively.

Figure 2(a) shows the comparison between the numerical results and
Eq.(\ref{psz2}) for zeros located on a window around $E=2.5041178 \times
10^{15}$ (as in Fig.~\ref{Fig1}). The effective matrix size is $N_{\rm eff} =
7.7376$ (instead of $N_{0} = 33.6188$), and $\alpha = 1.0438$. The agreement
is quite good, and shows that $N_{\rm eff}$ is the correct matrix size in this
case. For comparison, we have plotted as a dashed curve the theoretical
formula (\ref{psz2}) without the rescaling of the variable $s$.

\begin{figure}
\begin{center}
\leavevmode
\epsfysize=2.5in
\epsfbox{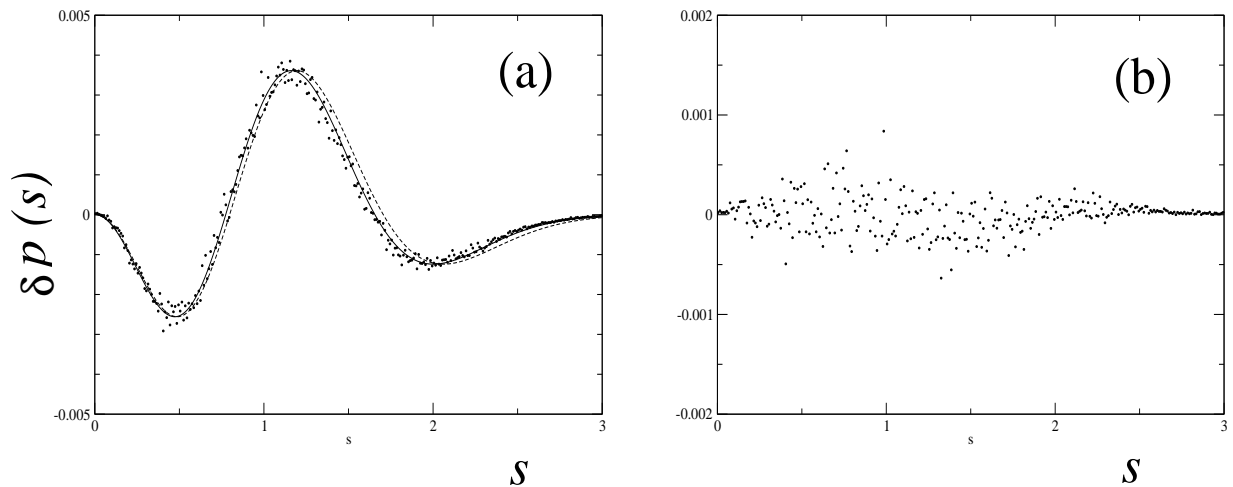}

\caption{{\small Same as in Fig.~2 but for a billion zeros located in a window
near $E=1.30664344 \times 10^{22}$.}}

\end{center}
\label{Fig3}
\end{figure}

Figure 2(b) is a plot of the difference between Odlyzko's results and the
prediction (\ref{psz2}). There is still some structure visible, which might be
attributed to the ${\cal O} ( N_{\rm eff}^{-4} )$ correction. To test the
convergence, we have made the same plot but now using one billion zeros
located on a window around $E=1.30664344 \times 10^{22}$, which corresponds to
$N_{\rm eff} = 11.2976$ (instead of $N_{0} = 49.0864$) and $\alpha = 1.0300$
(Fig.~\ref{Fig3}(a)). Now the agreement is clearly improved. The difference
between the prediction (\ref{psz2}) and the numerical results, plotted in
Fig.~\ref{Fig3}(b), shows a structureless remain.

\section{Conclusion}

For zeros of the Riemann zeta function located around $z = 1/2 + \ui E$, we
derived an heuristic formula for the nearest--neighbour spacing distribution
that contains finite--$E$ corrections. We argued that to leading order the
corrections are the same as for random matrices from the circular unitary
ensemble of size $N_{\rm eff} \approx 0.230158 \log (E/2\pi)$. We proposed to
describe the next to leading order correction as a simple scaling of the
dominant term.

Two main conjectures were used. The most important is the explicit expression
for the two--point correlation function for the Riemann zeros
(Eqs.~(\ref{diag}) and (\ref{off})) obtained in \cite{bk}. The second is the
statement that to leading order deviations from random matrix predictions
reduce to a change of the kernel (\ref{kern}). It follows that our procedure
to compute the finite--$E$ corrections of the complex Riemann zeros from
unitary matrices should be valid for any local statistics. 

It is known that at finite heights $E$ above the real axis the appropriate
symmetry group of the Riemann zeta function is not ${\rm CUE_N}$. It has been
conjectured \cite{ksa} that $\zeta (z)$ is a member of a family obeying the
symplectic symmetry $U Sp (2N)$. The possible connections of our results with
finite--$E$ corrections associated to this symmetry deserve to be
investigated.

We acknowledge discussions and data from A. Odlyzko.

\end{document}